# In memoriam Pál Révész (1934-2022).

## The random walk from orthogonality to anisotropy.

**Endre Csáki**
Alfréd Rényi Institute of Mathematics, Budapest, P.O.B. 127, H-1364, Hungary. E-mail address: csaki.endre@renyi.hu

**Antónia Földes**
Department of Mathematics, College of Staten Island, CUNY, 2800 Victory Blvd., Staten Island, New York 10314, U.S.A. E-mail address: antonia.foldes@csi.cuny.edu

**Abstract**
A short account of Pál Révész wonderful achievements in probability.

## 1 Introduction

The Hungarian mathematical community lost one of his leading members, when Pál Révész passed away on 14 of November 2022. He studied mathematics at the Eötvös Loránd University and graduated in 1957. He was a student and later a friend of Alfréd Rényi. His first appointment was in his alma mater as an assistant professor. He exchanged his teaching job in 1963 for a position at the Mathematical Research Institute (now Rényi Institute), where he worked as a research fellow and later as the head of the Department of Probability until 1985. Then he became a professor and Head of the Department of Statistics and Probability of the Vienna University of Technology, from where he retired in 1998 and was deputy chairman till 2005. From Vienna he regularly visited the Rényi Institute to talk mathematics with his friends and former colleagues. In all these years he traveled and taught in many places: in Addis Ababa, Ethiopia, at the University of Leiden in the Netherland, the ETH Zurich in Switzerland. He was a honorary professor at Carleton University in Ottawa Canada and the University of Szeged, where he taught regularly. And of course he was invited and gave talks all over Europe and the United States. He became a member of the Hungarian Academy 1982. His honors and achievements are hard to list. Just to name a few: he was president of the Bernoulli Society for Mathematical Statistics and Probability from 1983 to 1985, he became a member of the Academy Europaea in 1991, got the State Award of the People's republic of Hungary. He wrote four books, the third one: Random walk in random and non-random environments, was so successful, that it had three editions (1990, 2005 and 2013).

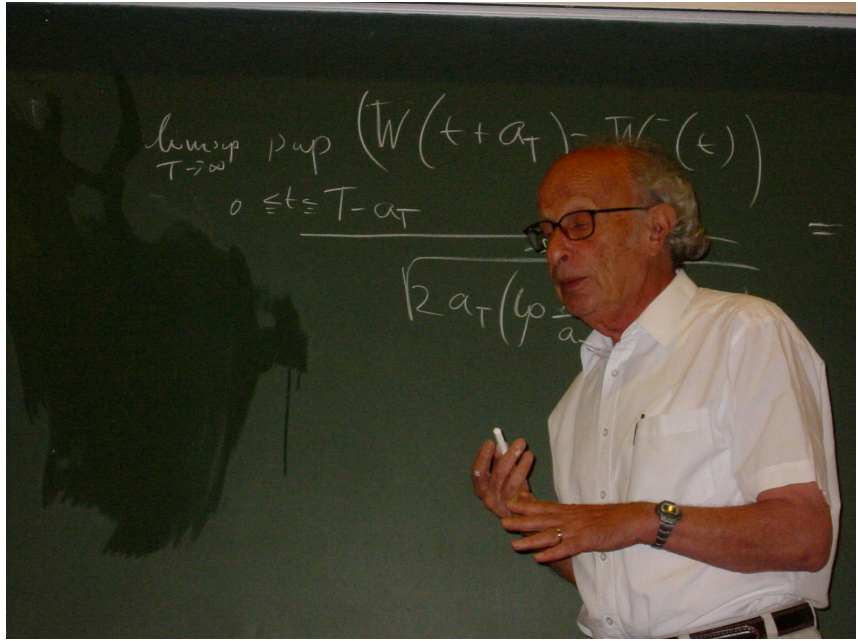

Figure 1: Budapest 2007

There were four conferences celebrating his 65th birthday in Balatonlelle in 1999, and three other ones in Budapest celebrating the 70th and the 75th and 80th birthday of Pál and Endre, the first author of this paper.

On each of these occasions Miklós Csörgő, Pál's friend and collaborator for almost half a century, gave a lecture. In those lectures, Miklós collected Pál most important results; [11]- [14]. In addition the second author of this paper, Antónia also gave a celebratory lecture at each occasion [17]-[19]. Miklós and the authors of this paper were collaborating so about half a century. Even though, for the last 30 years we all have been living in different countries, we spent a big part of the summers together in Budapest and discussed our research plan for the year ahead, and of course we did the rest by continuous emailing to each other. There are approximately 70 papers, which was written by Pál and a subset of the three of us. So here we will only give a few glimpses of Paul's work, mostly quoting from Miklós, with his permission, and his encouragement and present a more extensive picture of those problems, on which we worked together in the last ten plus years (thus not collected in Miklos's celebratory papers).

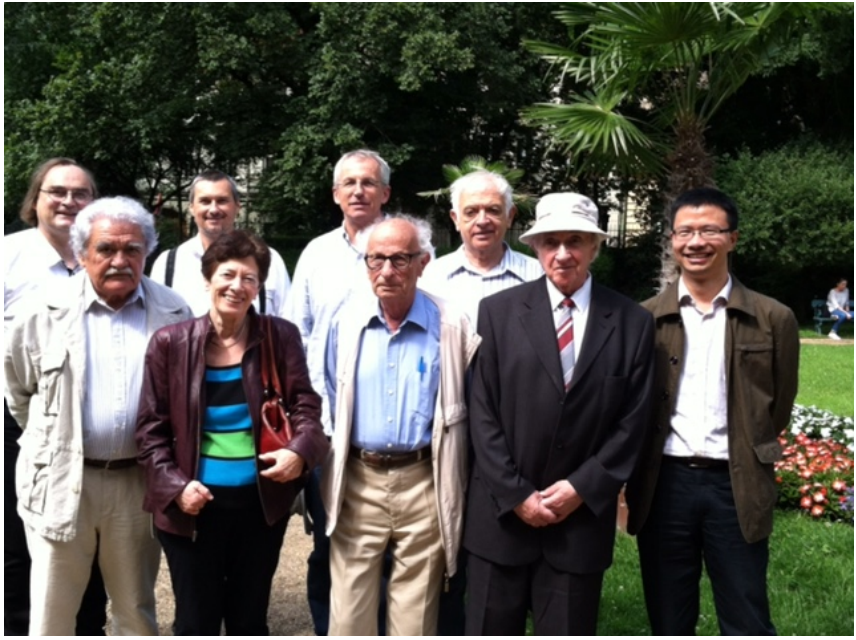

Figure 2: Budapest 2013

## 2   A few glimpses into the work of Pál Révész

Our choice of topics here, are the ones we treasure most, without trying to paint the full picture, which was done by Miklós. Let's start with Pál's first book; *The Laws of Large Numbers*, published in 1968. This book contains many important results about symmetrically dependent and orthonormal systems of random variables. From the point of view of a probabilist, these type of random variables are interesting as they have similar properties to the independent random variables. One of the central questions was the generalization of the well-known Kolmogorov three series theorem for independent random variables. This was achieved for orthonormal random variables by the Rademacher- Mensov theorem. Alexits and Tandori (1961) proposed to study the equinormed strongly multiplicative system (ESMS) of random variables: $\xi_1, \xi_2, \ldots$ is an ESMS if

$$\mathbf{E}\left(\xi_i\right) = 0, \quad \mathbf{E}\left(\xi_i^2\right) = 1 \quad i = 1, 2, \ldots \quad \mathbf{E}\left(\xi_{i_1}^{r_1} \xi_{i_2}^{r_2} \cdots \xi_{i_k}^{r_k}\right) = \mathbf{E}\left(\xi_{i_1}^{r_1}\right) \mathbf{E}\left(\xi_{i_2}^{r_2}\right) \cdots \mathbf{E}\left(\xi_{i_k}^{r_k}\right) \quad k = 1, 2, \ldots$$

where $r_1, r_2 \cdots r_k$ can be equal to 1 or 2. They proved that:

*For a uniformly bounded sequence of $\{\xi_k, \quad k = 1, 2 \ldots\}$ ESMS random variables*

$$\sum_{k=1}^{\infty} c_k \xi_k < \infty \quad \text{with probability 1.}$$

*with any sequence of real numbers $c_1, c_2, \ldots$ for which $\sum_{k=1}^{\infty} c_k^2 < \infty$.*

Révész proved the first law of iterated logarithm (LIL) and the first central limit theorem (CLT) for uniformly bounded ESMS system of random variables. Somewhat later, studying multiplicative systems in the seventies, he proved the first complete LIL for a sequence of uniformly bounded random variables $\{\xi_i, i = 1, 2, \ldots\}$ for which

$$\mathbf{E} \prod_{i=1}^{j} \xi_{\ell(i)} = 0 \quad \text{and} \quad \mathbf{E} \prod_{i=1}^{j} \xi_{\ell(i)}^2 = 1 \quad \text{for each} \quad j \geq 1$$

where $\ell(1) < \ell(2) < \cdots < \ldots \ell(j)$. Namely he proved that:

*Under the above conditions*

$$\limsup_{n \to \infty} \frac{\sum_{i=1}^{n} \xi_i}{(2n \log \log n)^{1/2}} = 1 \quad \text{a.s.}$$

One of the highlights of Pál's early work is in his paper "A problem of Steinhaus " [R21]. The glorious history of this problem is nicely described by Miklós in [12]. In a nutshell, Steinhaus conjecture turned out to be untrue, the first few steps were made by D.G. Austin, D. Darling and Rényi, but the best resolution was given by Pál. Instead of giving his general theorem we only quote his beautiful corollary:

*From any orthogonal sequence $f_n \in L^2$ one can always choose a subsequence $\{f_{n_k}\}$ such that $\sum_{k=1}^{\infty} c_k f_{n_k}$ converges a.e. and $L^2$ whenever $\sum_{k=1}^{\infty} c_k^2 < \infty$.*

In the seventies Pál turned to some problems of practical interest. The two most important issues from these topics were the density estimation and the Robbins -Monro type stochastic approximation. He wrote about twenty papers on these topics. As to the density estimation, the most representative example is a general definition of the empirical density. Let $X_1, X_2, \ldots$ be an independent identically distributed (i.i.d.) sequence of random variables with density function $f(x)$. Let $f(x)$ be vanishing outside $-\infty \leq C < D \leq \infty$ and let $\psi_k(x, y), \quad k = 1, 2, \ldots$ be a sequence of Borel measurable functions defined on the square $(A, B)^2$ with $-\infty \leq A \leq C < D \leq B \leq \infty$. Now a possible definition of the empirical density function is

$$f_n(x) = \frac{1}{n} \sum_{k=1}^{n} \psi_n(x, X_k) = \int_C^D \psi_n(x, X_k) \, dF_n(y)$$

$$F_n(y) = \frac{1}{n} \sum_{k=1}^{n} \mathbf{1}_{(-\infty, y]}(X_i)$$

where $F_n(y)$ is of course the empirical distribution function, based on the first $n$ observations of $F(y) = \int_{-\infty}^{y} f(x) dx$. This definition as it is explained in details in the paper [R42] (the first joint work of Pál and Antónia) contains many earlier type of density estimations, the histogram, the kernel type (Rosenblatt and Parzen), the orthogonal expansion type definition, cf. Čencov, van

Ryzin, and Schwartz. In paper [R42] conditions are given for the nearness in sup norm of $f_n(x)$ to $f(x)$ at an exponential rate in probability.

Révész had a huge contribution to the theory of Robbins-Monro type stochastic approximation. See [R34], [R37], [R41], [R44], [R45], [R46], [R57] and [R58]. Namely he was the first to get an estimation of regression function simultaneously in all its points. In [R45] he starts with $(X_1, Y_1), (X_2, Y_2), \ldots$ i.i.d. pairs of random variables such that $X \in [0, 1], Y \in R^1$ and then considers the regression function $r(x) = \mathbf{E}(Y|X = x)$. Starting with an arbitrary continuous function $r_0(x)$, Pál defined his Robbins-Monro type recursive estimation of the regression function as

$$r_{n+1}(x) = r_n(x) + \frac{1}{(n+1)a_{n+1}} K\left(\frac{x - X_{n+1}}{a_{n+1}}\right)(Y_{n+1} - r_{n+1}) \quad 0 \leq x \leq 1.$$

where $a_n = n^{-\alpha}$, $0 < \alpha < 1$, and $K(.)$ is a density function. Under some natural condition he gives an exponential bound for the nearness of $r_n(x)$ to $r(x)$ in $L^2$, proves a pointwise and an $L^2$ CLT (central limit theorem) for $r_n(x) - r(x)$ and a.s. uniform convergence with rates:

$$\mathbf{P}\left(\sup_x |r_n(x) - r(x)| \geq \frac{C \log n}{(na_n)^{1/2}}\right) = O\left(n^{-T}\right)$$

*where $C = C(T) > 0$ is a constant depending on $T$ only.*

Révész collaboration with Pál Erdős started around 1972. They wrote at this time two joint papers [R56a] and [R56b] about the length of the longest pure head run (uninterrupted sequence of heads) in the coin tossing game with a fair coin. To give a taste of their beautiful results we have to recall some definitions.

The function $a_1(t)$ belongs to the upper-upper class of the stochastic process $\{Y(t)\}$ abbreviated as $(a_1(t) \in UUC(Y(t)))$ if for almost all $\omega \in \Omega$ there exists a $t_0(\omega) > 0$ such that $Y(t) < a_1(t)$ if $t > t_0 = t_0(\omega)$.

The function $a_2(t)$ belongs to the upper-lower class of the stochastic process $\{Y(t)\}$ abbreviated as $(a_2(t) \in \mathrm{ULC}(Y(t)))$ if for almost all $\omega \in \Omega$ there exists a sequence of positive numbers $0 < t_1 = t_1(\omega) < t_2 = t_2(\omega) \ldots$ with $t_n \to \infty$ such that $Y(t_i) \geq a_2(t_i)$ $(i = 1, 2, \ldots)$.

Similarly we get the definition of the lower- upper class $a_3(t) \in LUC$ by repeating the second definition and changing the inequalities into $Y(t_i) \leq a_3(t_i)$ $(i = 1, 2, \ldots)$.

Finally to get the definition of the lower- lower class $a_4(t) \in LLC$ we need to repeat the first definition and modify the inequalities into $Y(t) > a_4(t)$ if $t > t_0 = t_0(\omega)$.

These definitions can be used for sequences, instead of function as as well. In fact this is what see in the following theorem. Erdős and Révész proved the following amazingly precise result:

*Let $Z_n$ be the length of the longest pure head run in the first $n$ tosses of a fair coin, and let*

$$B(\{b_n\}) = \sum_{n=1}^{\infty} 2^{-b_n},$$

*Then*
$$b_n \in UUC(Z_n) \quad \text{if} \quad B(\{b_n\}) < \infty$$

$$b_n \in ULC(Z_n) \quad \text{if} \quad B(\{b_n\}) = \infty$$

*and for any $\epsilon > 0$*

$$\kappa_n := [\log_2 n - \log_2 \log_2 \log_2 n + \log_2 \log_2 e - 1 + \epsilon] \in LUC(Z_n),$$
$$\lambda_n := [\log_2 n - \log_2 \log_2 \log_2 n + \log_2 \log_2 e - 2 + \epsilon] \in LLC(Z_n).$$

Here the first two statements are the best possible results. The third and fourth, the lower class results are nearly best possible, their complete characterization was given by Guibas-Odlyzko [20] and Samarova [40].

Pál and Miklós collaboration started in 1975 with a problem of Strassen [44]. Strassen showed that:

*On the same probability space a sequence of i.i.d. random variables $X_1, X_2, \ldots$, with $\mathbf{E}(X_1) = 0$, $\mathbf{E}(X_1^2) = 1$, $\mathbf{E}(|X_1|^4) < \infty$ and a standard Wiener process $\{W(t) \quad 0 < t \leq \infty\}$ can be constructed, such that as $n \to \infty$,*

$$Z_n := |S(n) - W(n)| = O\left((n \log \log n)^{1/4} (\log n)^{1/2}\right) \quad \text{a.s.},$$

*where $S(n) = \sum_{i=1}^n X_i$.*

He also posed the following question. Is it true that if we change $O(\cdot)$ into $o(\cdot)$, then it implies that $X_1$ is standard normal?

J. Kiefer [30] proved, that using the so called Skorohod [42] embedding scheme Strassen's above result cannot be improved. In their paper [R47] Pál and Miklós established the following strong approximation result:

*Assuming Cramer's condition on the characteristic function of the random variable $X_1$, for any given $0 < \varepsilon < 1/2$ one can assume enough moment conditions for $X_1$ of the i.i.d. sequence with $E(X_1) = 0$, $E(X_1^2) = 1$, $E(X_1^3) = 0$, .... such that with an appropriately constructed Wiener process $\{W(t) \quad 0 < t < \infty\}$ on the same probability space, as $n \to \infty$ one has*

$$|S(n) - W(n)| = o(n^\varepsilon) \quad \text{a.s.} \tag{2.1}$$

This was achieved with a new method of construction, based on the so called quantile transformation, which was first developed in [1] by Bártfai in 1966.

Kiefer [31] was the first one, who constructed an almost sure representation of the empirical process, by embedding the empirical process into an appropriate two-time parameter Gaussian process via his ingenious extension of the Skorohod embedding scheme to the case of of vector valued random variables. If $U_1, U_2, \ldots U_n, n = 1, 2, \ldots$ are independent random d-vectors, uniformly distributed on $[0,1]^d$, $d \geq 1$ one can define the uniform empirical process of these variables by

$$\alpha_n(y) := \sqrt{n}\left(\frac{\sum_{i=1}^n \mathbb{1}_{(0,y]}(U_i)}{n} - \lambda(y)\right), \quad y \in [0,1]^d,$$

where $\lambda(y)$ is the Lebesgue measure of the $d$-dimensional interval $(0, y]$. A separable mean zero Gaussian process defined on $[0,1]^d \times [0, \infty)$, $\{K(y,t); y \in [0,1]^d, t \in [0, \infty)\}$, is called a Kiefer process if

$$\mathbf{E}(K(x,s)K(y,s)) = (\lambda(x \wedge y) - \lambda(x)\lambda(y))(s \wedge t),$$

where the minimum $x \wedge y$ is taken component-wise and $\lambda(.)$ is as above, i.e. the Gaussian process $K(\cdot, \cdot)$ is a Brownian bridge in its first argument and a Wiener process in its second argument. In [31] Kiefer proved that:

*One can construct a probability space for $U_1, U_2, \ldots$ with a Kiefer process $K(\cdot, \cdot)$ on it such that as $n \to \infty$*

$$\max_{1 \leq k \leq n} \sup_{y \in [0,1]^d} |\sqrt{k}\alpha_k(y) - K(y,k)| = O(n^{1/3}(\log n)^{2/3}) \quad a.s.$$

In [R48] Pál and Miklós combined the classical Poissonization technique with their quantile transform method to prove the following important result:

*One can construct a probability space for $U_1, U_2, \ldots$ with a sequence of Brownian bridges $\{B_n(y); y \in [0,1]^d, d \geq 1\}_{n=1}^\infty$ on it, so that for any $r > 0$ there is a constant $C > 0$ such that for each $n = 2, 3 \ldots$*

$$\mathbb{P}\left(\sup_{y \in [0,1]^d} |\alpha_n(y) - B_n(y)| > C(\log n)^{3/2} n^{-1/2(d+1)}\right) \leq n^{-r}, \quad d \geq 1, \tag{2.2}$$

*and with a Kiefer process $\{K(y,t); (y,t) \in [0,1]^d \times [0, \infty)\}$, on it, so that for any $r > 0$ there is a constant $C > 0$ such that for all $n \geq 2$*

$$\mathbb{P}\left(\max_{1 \leq k \leq n} \sup_{y \in [0,1]^d} |\sqrt{k}\alpha_k(y) - K(y,k)| > Cn^{(d+1)/2(d+2)}(\log n)^2\right) \leq n^{-r}, \quad d \geq 1. \tag{2.3}$$

Actually the name Kiefer process was coined by Pál and Miklós in [R48] as well. The above two coupling type inequalities for the uniform empirical process and their approximating Gaussian counterpart are first of their kind in the literature. We do not want to delve into the many refinements of (2.1) - (2.3) of [R47] and [R48]. However the dyadic scheme refinement of their quantile transform method, which was developed by Komlós, Major, and Tusnády (1975, 1976) in [32] and [33], became now known as the famous KMT method or Hungarian construction.
Without quoting all the amazingly precise results of the KMT method, we only mention that they proved:

*it is enough to assume that $X_1$ has mean zero, variance one, and has $p > 2$ moments, that (2.1) should hold with $\varepsilon = 1/p$.*

This whole wonderful story is detailed in [12].

One of the highlights of Pál and Miklós collaboration is their famous *Strong Approximation in Probability and Statistics* book (1981). In this book Chapters 2-6 deal with this topic and its many generalizations as well. In their papers [R60] and [R63] they studied the fine analytic properties of the increments of the Wiener process and other related processes. These issues and its fascinating consequences for the sums of i.i.d. random variables using the strong approximation method are all presented in their book. We just quote one of their beautiful and incredible precise results, which proved to be extremely useful over the years:

*Let $a_T$ be a monotonically non-decreasing function of $T$ such that $0 < a_T \leq T$ and $T/a_T$ is monotonically non-decreasing as well. Define*

$$\beta_T := \left(2a_T \left(\log \frac{T}{a_T} + \log\log T\right)\right)^{-1/2}.$$

*Then*

$$\limsup_{T \to \infty} \sup_{0 < t \leq T - a_T} \beta_T |W(t + a_T) - W(t)| = 1 \quad \text{a.s.},$$

$$\limsup_{T \to \infty} \beta_T |W(T + a_T) - W(T)| = \limsup_{T \to \infty} \sup_{0 \leq s \leq a_T} \beta_T |W(T + s) - W(T)| = 1 \quad \text{a.s.},$$

*and*

$$\limsup_{T \to \infty} \sup_{0 < t \leq T - a_T} \sup_{0 \leq s \leq a_T} \beta_T |W(t + s) - W(t)| = 1 \quad \text{a.s.}$$

*If we also have*

$$\lim_{T \to \infty} \log(T/a_T) / \log\log T = \infty,$$

*then $\limsup_{T \to \infty}$ can be replaced by $\lim_{T \to \infty}$ in all the above statements.*

The above result for example contains the LIL for the Wiener process, and a special case of the Erdős- Rényi law, which states that:

*For any $c > 0$*

$$\lim_{n \to \infty} \max_{0 \leq k \leq n - [c \log n]} \frac{S(k + [c \log n]) - S(k)}{[c \log n]} = \alpha(c)$$

*where $S(n)$ is the sum of $X_1, X_2, \ldots X_n$ i.i.d. random variables, $\alpha(c) := \sup\{x : \rho(x) \geq e^{-1/c}\}$, and $\rho(x)$ is the Chernoff function of $X_1$. Furthermore the function $\alpha(c)$ uniquely determines the moment generating function and hence also the distribution function of $X_1$.*

However if the window $a_T$ is somewhat bigger, that is $a_T/T \to \infty$, then the increments of the partial sum $S(n)$ behave like the increments of the Brownian motion, by strong approximation and KMT, as long as $X_1$ has a moment generating function.

Going back to the increment problem above, the issue of the liminf of these increments brought together Pál and Endre for the first time. In [R69] they proved many interesting results, from which we only quote the following:

*Suppose that $0 < a_T \leq T$ and $a_T/T$ is non-increasing and*

$$\lim_{T \to \infty} \frac{\log (T/a_T)}{\log \log \log T} = \infty.$$

*Let*

$$\gamma(T) = \left( 2a_T \log \left( 1 + \frac{\pi^2}{16} \frac{[T/a_T]}{\log \log T} \right) \right)^{-1/2}.$$

*Then*

$$\liminf_{T \to \infty} \gamma(T) \sup_{0 < t \leq T - a_T} \sup_{0 \leq s \leq a_T} |W(t+s) - W(t)| = 1.$$

Let us mention one more very important result of Pál and Miklós about the small increments of the Wiener process, [R63]:

$$\lim_{h \to 0} \inf_{0 \leq s \leq 1-h} \sup_{0 \leq t \leq h} \sqrt{\frac{8 \log h^{-1}}{\pi^2 h}} |W(s+t) - W(s)| \stackrel{a.s.}{=} 1.$$

This theorem implies the well-known fact:

*Almost all sample functions of a Wiener process are nowhere differentiable.*

They also developed a lot of application of their strong invariance results. There are two chapters in their book, Chapter 4 and 5, see also their paper [R53], where they approximate the empirical processes by Gaussian processes and study the empirical and the quantile processes, using strong approximation methods. One of the highlights of this topic is Theorem 4.5.7 in their book.

In the end of the seventies Pál got interested in local time problems. So much so, that investigation of the properties of the local time is taking a central stage in his most successful third book; *Random walk in Random and Non-Random Environments*. Consider a simple symmetric walk, that is $S(0) = 0, S(k) = \sum_{i=1}^{k} X_i$, $k \geq 1$, where $\{X_i, i = 1, 2, \ldots\}$ are i.i.d. random variables with $\mathbf{P}(X_1 = 1) = \mathbf{P}(X_1 = -1) = 1/2$. The local time of the walk $S(.)$ is defined as

$$\xi(x,n) := \#\{k : 1 \leq n, S(k) = x\}, \quad x = 0, \pm 1, \pm 2, \ldots. \tag{2.4}$$

which is the number of times that the walk is in the position $x$ up to time $n$. The local time of the Brownian motion is defined as the Radon-Nikodym derivative $\{L(x,t), x \in R \quad t \geq 0\}$ for the occupation time as follows:

$$H(A,t) := \lambda\{s : 0 \leq s \leq t, W(s) \in A\} = \int_A L(x,t)dx$$

for any $t > 0$, and any Borel set $A$ of the real line, where $\lambda(.)$ is the Lebesgue measure. Trotter proved in 1958 [45] that $L(x,t)$ is continuous in both arguments. He also proved limsup results for the small increments of $L(x,t)$ in $x$ and the $t$ variable as well. These were improved by H.P. Mc Kean Jr. [35], D.B. Ray [39], J. Hawkes [22] and E. Perkins [37]. H. Kesten [28] proved the following LIL:

$$\limsup_{t\to\infty} \frac{L(0,t)}{(2t\log\log t)^{1/2}} = \limsup_{t\to\infty} \frac{\sup_{-\infty<x<\infty} L(x,t)}{(2t\log\log t)^{1/2}} = 1 \quad \text{a.s.}$$

The first paper the four of us wrote together is [R88] in Pál's list, was about the big increments of the local time in 1983. Using the notation $L(t) = L(0,t)$ we proved an analogous result to the Wiener increment theorem above, namely:

Let $0 < a_t \leq t$ be a non-decreasing function of $t$. Assume that $a_t/t$ is non-increasing. Then

$$\limsup_{t\to\infty} \gamma_t Y(t) = \limsup_{t\to\infty} \gamma_t a_t^{-1/2}\left(L(t) - L\left(t - a_t\right)\right) = 1 \quad \text{a.s.}$$

*where*

$$Y(t) = Y(t, a_t) = a_t^{-1/2} \sup_{0 \leq s \leq t - a_t} \left(L\left(s + a_t\right) - L(s)\right)$$

*and*

$$\gamma_t = \left(\log\left(t/a_t\right) + 2\log\log t\right)^{-1/2}.$$

*If we also assume that*

$$\lim_{t\to\infty} \frac{\log\left(t/a_t\right)}{\log\log t} = \infty,$$

*then*

$$\lim_{t\to\infty} \gamma_t Y(t) = 1 \quad \text{a.s.}$$

Pál was the first, who proved a strong invariance principle for the local time with rate, using the Skorohod construction. He established in [R77] that:

*On a rich enough probability space for any $\varepsilon > 0$ as $n \to \infty$*

$$\sup_x |\xi(x,n) - L(x,n)| = o\left(n^{1/4+\varepsilon}\right),$$

*where the sup is taken for integers.*

It turned out that this result is very close to the best possible one. Namely, the rate $o\left(n^{1/4+\varepsilon}\right)$, cannot be replaced by $o\left(n^{1/4}\right)$. Dobrushin [25] proved in 1955 that:

$$\frac{\xi(1,n) - \xi(0,n)}{2^{1/2} n^{1/4}} \xrightarrow{d} Z_1 |Z_2|^{1/2},$$

*and for the Brownian local time we have* (see e.g. Yor [48] )

$$\frac{L(1,n) - L(0,n)}{2 n^{1/4}} \xrightarrow{d} Z_1 |Z_2|^{1/2},$$

*where $Z_1$ and $Z_2$ are independent standard normal random variables.*

That is to say we have the same limit but with a little different normalization. However $o\left(n^{1/4+\varepsilon}\right)$ can be replaced by $O\left(n^{1/4} \log n\right)$ as it was shown by Borodin [6]. One of the most useful and beautiful results of Pál is, that he managed to prove a joint approximation result for random walk and the Wiener process and their respective local times as well [R77]:

*On an appropriate probability space for a simple symmetric random walk $\{S(n); n = 1, 2, \ldots\}$ with local time $\{\xi(x,n); x = 0, \pm 1, \pm 2 \ldots\}$ one can construct a standard Wiener process $\{W(t); t \geq 0)\}$ with local time process $\{\eta(x,t) : x \in R, t \geq 0\}$ such that, as $n \to \infty$ we have for any $\epsilon > 0$*

$$|S(n) - W(n)| = O\left(n^{1/4+\varepsilon}\right) \quad \text{a.s.}$$

*and*

$$\sup_{x \in Z} |\xi(x,n) - \eta(x,n)| = O\left(n^{1/4+\varepsilon}\right) \quad \text{a.s.}$$

*hold simultaneously.* Just to mention a few more topics from this book, Pál investigates the local time of the random walk and the Wiener process in one dimension and $\mathbb{Z}^n$, the exact and limit distributions for them, the connection of the local time and the high excursions, and the strong approximation of the local time difference. Let us present the following example [R116] which is actually a joint work of the four of us: *There exists a probability space with*

- *a standard Wiener process $\{W(t); t \geq 0\}$ and its two parameter local time process $\{\eta(x,t); x \in \mathbb{R}^1, t \geq 0\}$*

- *a two-time parameter Wiener process $\{W(x,u); x \geq 0, u \geq 0\}$*

- *a process $\{\hat{\eta}(0,t); t \geq 0\} \stackrel{\mathbb{D}}{=} \{\eta(0,t); t \geq 0\}$*

*such that*

$$\sup_{0 \leq x \leq A t^{\delta/2}} |\eta(x,t) - \eta(0,t) - 2W(x, \hat{\eta}(0,t))| = O\left(t^{(1+\delta)/4 - \varepsilon/2}\right) \quad \text{a.s.} \quad (t \to \infty), \qquad (2.5)$$

$$\{\hat{\eta}(0,t); t \geq 0\} \quad \text{and} \quad \{W(x,u); x \geq 0, u \geq 0\} \quad \text{are independent},$$

where

$$A > 0, \quad 0 \leq \delta < 7/100, \quad 0 < \varepsilon \leq 1/72 - \delta/7.$$

The strength of this result is the independence of the processes $\{W(x,u); x \geq 0, u \geq 0\}$ and $\{\hat{\eta}(0,t); t \geq 0\}$ and the fact that we know a lot about these processes. The process $W(x, \hat{\eta}(0,t))$ is called an iterated Wiener process, which has many nice properties. In the book there is one and a half page containing just the consequences of the above theorem. Here we mention only the following pair:

$$\limsup_{t \to \infty} \frac{\eta(x,t) - \eta(0,t)}{2\sqrt{2x\eta(0,t) \log \log t}} = 1 \quad \text{a.s.} \quad \text{for any} \quad x \geq 0,$$

$$\limsup_{t \to \infty} \frac{\eta(x,t) - \eta(0,t)}{x^{1/2} t^{1/4} (\log \log t)^{3/4}} = \frac{4}{3} 6^{1/4} \quad \text{a.s.} \quad \text{for any} \quad x \geq 0.$$

Another beautiful area of Pál's investigation is the topic of excursions. In fact he has many results which investigate the strong connection between local time and excursions. Consider a simple symmetric random walk. Let $\rho_0 = 0$,

$$\rho_1 = \inf\{k : k > 0 \ S(k) = 0\}, \quad \rho_2 = \inf\{k : k > \rho_1, S(k) = 0\}, \ldots \rho_{n+1} = \inf\{k : k > \rho_n, S(k) = 0\}.$$

We call the parts of the walk $\{S(\rho_{i-1}+1), S(\rho_{i-1}+2), \ldots, S(\rho_i)\}$ $i = 1, 2, \ldots$ excursions. So $\rho_{i+1} - \rho_i$ is the length of the i-th excursion. Pál and Endre were interested in the length of the longest excursion. They wrote a paper [R96a] for the 70-th birthday of Paul Erdős. We only quote the following interesting result. Let

$$R(n) = \#\{k : 0 < k \leq n, S(k) > 0\}, \quad T(n) = \max(\rho_1, \rho_2 - \rho_1, \ldots, \rho_{R(n)-1} - \rho_{R(n)}, n - \rho_{R(n)})$$

$R(n)$ is the number of excursions up to $n$, and $T(n)$ is the length of the longest excursion, including the last unfinished one up to $n$. They proved that: *For the longest excursion $T(n)$*

$$\liminf_{n \to \infty} \frac{\log \log n}{n} T(n) = \beta \quad \text{a.s.}$$

*where $\beta$ is the only solution of the equation of*

$$\sum_{k=1}^{\infty} \frac{\beta^k}{k!(2k-1)} = 1. \tag{2.6}$$

Consider now the excursion lengths $\rho_1, \rho_2 - \rho_1, \ldots, \rho_{R(n)-1} - \rho_{R(n)}, n - \rho_{R(n)}$ and create its ordered sample

$$T(n) = T_1(n) \geq T_2(n) \geq \ldots \geq T_{R(n)+1}(n)$$

Pál and Endre together with Erdős showed in [R96b] that: *For any fixed k=1,2, ... we have*

$$\liminf_{n\to\infty} \frac{\log\log n}{n} \sum_{j=1}^{k} T_j(n) = k\beta \quad a.s.$$

There is a whole chapter, titled Excursion, with detailed description of these two papers and many related issues in Pál's third book. Let us mention one more thing, which connects the above results on excursions with the earlier mentioned iterated processes. Recall the definition of $\xi(0,n)$, the local time of zero in (2.4). Consider now the random walk counterpart of the iterated process in (2.5). Let $\{S_1(n), n = 1, 2, \ldots\}$ a simple symmetric random walk and let $\xi_2(k,n)$ $k = 0, \pm 1, \pm 2, \ldots$ $n = 1, 2, \ldots\}$ the local time process of a second simple symmetric random walk which is independent from $S_1(.)$, and define
$R(n) = S_1(\xi_2(n-1))$ $n = 1, 2, \ldots$ and the occupation time of $R(n)$ as

$$\xi^*(r,n) = \#\{k : 1 \le k \le n, R(k) = r\}, \ r = 0, \pm 1, \pm 2, \ldots$$

In [R151] the four of us investigated the properties of $\xi^*(r,n)$ and its Brownian counterpart as well. Here we just quote the following pair of results:

*We have*
$$\limsup_{n\to\infty} \frac{1}{n} \max_r \xi^*(r,n) = 1 \quad a.s.$$

*and*
$$\liminf_{n\to\infty} \frac{\log\log n}{n} \max_r \xi^*(r,n) \ge \beta \quad a.s. \tag{2.7}$$

*where $\beta$ is the same as in (2.6).*

We also stated as a conjecture, that (2.7) is true with equality.

The random walk in random scenery model goes back to Kesten and Spitzer [29]. Here is their definition. Let $\mathbf{S}_n$ be a simple symmetric walk on $\mathbf{Z}^d$ starting from the origin. Whenever this random walk visits a site $\mathbf{x} \in \mathbf{Z}^d$ the fortune changes by an amount of $Y(x) \in R$, that is to say the total value of $\mathbf{S}_n$ at step $n$ is

$$Z(n) := \sum_{j=0}^{n} Y(\mathbf{S}_j).$$

Supposing that $\{Y(\mathbf{x}), \mathbf{x} \in \mathbf{Z}^d\}$ is a collection of i.i.d.r.v.-s with zero expectation and finite variance $\sigma^2$, this collection is called a random scenery which is supposed to be independent from the walk $\{\mathbf{S}_n\ n \ge 0\}$, while $\{Z(n)\ n \ge 0\}$ is called a random walk in a random scenery. As to the history of this topic we refer to [13]. Pál contributed to this work together with Zhan Shi in [R160] and together with Endre the three of them in [R165] proved among many other results the following strong approximation theorem:

*On a rich enough probability space, as $n \to \infty$*

$$\left| Z(n) - \sigma \left(\frac{2}{\pi}\right)^{1/2} W(n\log n) \right| = o\left(n^{1/2}(\log n)^{3/8+\epsilon}\right) \quad a.s.$$

*for any $\varepsilon > 0$ as $n \to \infty$.*

In [R175] Jay Rosen, Zhan Shi, Pál and us studied the occupation measure of various sets for symmetric transient random walk in $Z^d$ with finite variances. Let $X_i$, $i = 1, 2, \ldots$ be a sequence of independent identically distributed (i.i.d.) random variables and $S_n = \sum_{i=1}^n X_i$. Let $\mu_n^S(A) = \sum_{j=0}^n \mathbf{1}_A(S_j)$ for all sets $A \subseteq Z^d$, denote the occupation measure of the set $A$, where $S_j$ is a symmetric transient random walk in $\mathbf{Z}^d$ with $d \geq 3$. Assuming that $S_n$ is not supported on any subgroup strictly smaller than $Z^d$, let $G(x) = P(S_n = x)$ the Green function for $\{S_n\}$. Finally for any finite $A \subseteq Z^d$ let $\Lambda_A$ denote the largest eigenvalue of the $|A| \times |A|$ matrix $G_A(x,y) = G(x-y)$, $x, y \in A$. We proved that:

*If $X_1$ has finite second moment then*

$$\lim_{n \to \infty} \sup_{x \in Z^d} \frac{\mu_n^S(x+A)}{\log n} = \lim_{n \to \infty} \sup_{0 \leq m \leq n} \frac{\mu_n^S(S_m + A)}{\log n} = -1/\log(1 - 1/\Lambda_A) \quad \text{a.s.}$$

In [R175] many other results and examples were also given.

Around 2004 Pál got interested in the question of the number of zeros of a two-parameter random walk. First he, Davar Khoshnevisan and Zhan Shi [R174] investigated the following question. Let $\{X_{i,j}\}_{i,j=1}^\infty$ denote i.i.d. random variables, taking the values $\pm 1$ with respective probabilities $1/2$, and consider the two-parameter random walk $\mathbf{S} := \{S(n,m)\}_{n,m,\geq 1}$ defined by

$$S(n,m) = \sum_{i=1}^n \sum_{j=1}^m X_{i,j} \quad \text{for } n, m \geq 1.$$

A lattice point $(i, j)$ is said to be a vertical crossing for the random walk $\mathbf{S}$, if $S(i,j)S(i,j+1) \leq 0$ Let $\mathbf{Z}_n$ denote the total number of vertical crossings in the box $[1, N]^2 \cap \mathbf{Z}^2$. They proved that:

*With probability one*

$$Z(N) = N^{3/2+o(1)} \quad \text{as} \quad N \to \infty.$$

In [R192] Pál and Khoshnevisan returned to this question, but they described the asymptotic behavior of two other "contour plotting algorithms". They considered the total number of zeros in $[1, N]^2 \cap \mathbf{Z}^2$, defined by

$$\gamma_N := \sum \sum_{(i,j) \in [0,N]^2} \mathbf{1}_{\{S(i,j)=0\}}$$

and the total number of on-diagonal zeros in $[1, 2N]^2 \cap \mathbf{Z}^2$,

$$\delta_N := \sum_{i=1}^N \mathbf{1}_{\{S(2i,2i)=0\}}.$$

Their main result is that:

*With probability one*

$$\gamma_N = N^{1+o(1)} \quad as \quad N \to \infty \quad and \quad \lim_{N \to \infty} \frac{\delta_N}{\log N} = \frac{1}{(2\pi)^{1/2}}.$$

## 3  Our joint work on anisotropic walks.

After spending so much time studying the random walk and Brownian motion, the four of us begin to explore the anisotropic random walk on $\mathbb{Z}^2$. $\{\mathbf{C}(N) = (C_1(N), C_2(N)) \, ; \, N = 0, 1, 2, \ldots\}$ on $\mathbb{Z}^2$ is an anisotropic random walk, if it has the following transition probabilities

$$\mathbf{P}(\mathbf{C}(N+1) = (k+1, j)|\mathbf{C}(N) = (k, j)) = \mathbf{P}(\mathbf{C}(N+1) = (k-1, j)|\mathbf{C}(N) = (k, j)) = \frac{1}{2} - p_j,$$

$$\mathbf{P}(\mathbf{C}(N+1) = (k, j+1)|\mathbf{C}(N) = (k, j)) = \mathbf{P}(\mathbf{C}(N+1) = (k, j-1)|\mathbf{C}(N) = (k, j)) = p_j,$$

for $(k, j) \in \mathbb{Z}^2$, $N = 0, 1, 2, \ldots$ We assume that $0 < p_j \leq 1/2$ and $\min_{j \in \mathbb{Z}} p_j < 1/2$, and $\mathbf{C}(0) = (0, 0)$. This is a generalization of the simple symmetric walk on the plane, where $p_j = 1/4$ $j = 0, \pm 1, \pm 2, \ldots$ The waste literature of this model would be hard to list. Moreover it also contains the comb model, which was our first topic of investigation, before considering the general anisotropic walk.

**Comb model**

The special anisotropic model, where $p_0 = 1/4$, $p_j = 1/2$, $j = \pm 1, \pm 2, \ldots$ in which case all horizontal lines are missing except the $x$-axis is called the two-dimensional comb. As to the literature of the comb model we refer to Weiss and Havlin [47], Bertacchi and Zucca [5], Bertacchi [2], Csáki *et al.* [R187]. Our results on the comb model [R187], [R193] were summarized by Antónia [19] at the occasion of the 75th birthday conference so we only give a short account. We begin with another more straightforward (but equivalent) definition: The 2-dimensional comb lattice $\mathbb{C}^2$, is obtained from $\mathbb{Z}^2$ by removing all horizontal edges off the $x$-axis. Then the walk $\mathbf{C}(n) = (C_1(n), C_2(n))$ is defined by the transition probability

$$p(\mathbf{u}, \mathbf{v}) := \mathbf{P}(\mathbf{C}(n+1) = \mathbf{v} \mid \mathbf{C}(n) = \mathbf{u}) = \frac{1}{\deg(\mathbf{u})}$$

for locations $\mathbf{u}$ and $\mathbf{v}$ that are neighbors on $\mathbb{C}^2$, where $\deg(\mathbf{u})$ is the number of neighbors of $\mathbf{u}$, otherwise $p(\mathbf{u}, \mathbf{v}) := 0$. Consequently, the non-zero transition probabilities are equal to $1/4$ if $\mathbf{u}$ is on the horizontal axis and they are equal to $1/2$ otherwise. Lets start with a theorem of Bertacchi. Defining the continuous time process $\mathbf{C}(nt) = (C_1(nt), C_2(nt))$ by linear interpolation from $(C_1(n), C_2(n))$, she proved for the comb model that

$$\left(\frac{C_1(nt)}{n^{1/4}}, \frac{C_2(nt)}{n^{1/2}}; t \geq 0\right) \xrightarrow{\text{Law}} (W_1(\eta_2(0, t)), W_2(t); t \geq 0), \quad n \to \infty,$$

where $W_1$, $W_2$ are two independent Wiener processes (Brownian motions) and $\eta_2(0,t)$ is the local time process of $W_2$ at zero, and $\xrightarrow{\text{Law}}$ denotes weak convergence on $C([0,\infty),\mathbb{R}^2)$ endowed with the topology of uniform convergence on compact intervals.

The limit process of the first component above is a so called iterated process for which

$$\frac{W_1(\eta_2(0,t))}{t^{1/4}} \stackrel{\text{Law}}{=} X|Y|^{1/2}, \quad t > 0 \text{ fixed,}$$

where $X$ and $Y$ are independent standard normal random variables.

We proved an almost sure version of her theorem:

*On an appropriate probability space for the random walk $\{\mathbf{C}(n) = (C_1(n), C_2(n)); n = 0, 1, 2, \ldots\}$ on $\mathbb{C}^2$, one can construct two independent standard Wiener processes $\{W_1(t); t \geq 0\}$, $\{W_2(t); t \geq 0\}$ so that, as $n \to \infty$, we have with any $\varepsilon > 0$*

$$n^{-1/4}|C_1(n) - W_1(\eta_2(0,n))| + n^{-1/2}|C_2(n) - W_2(n)| = O(n^{-1/8+\varepsilon}) \quad a.s.,$$

*where $\eta_2(0,\cdot)$ is the local time process at zero of $W_2(\cdot)$.*

We also proved a Strassen type result for the joint behavior of the limiting processes in the above theorem. To present it we have to recall some definitions. Let $\mathcal{S}$ be the Strassen class of functions, i.e., $\mathcal{S} \subset C([0,1], \mathbb{R})$ is the class of absolutely continuous functions (with respect to the Lebesgue measure) on $[0,1]$ for which

$$f(0) = 0 \quad \text{and} \quad \int_0^1 \dot{f}^2(x)dx \leq 1.$$

and let $\mathcal{S}_M \subset \mathcal{S}$, the class of non-decreasing functions in $\mathcal{S}$.

Define the Strassen class $\mathcal{S}^2$ as the set of $\mathbb{R}^2$ valued, absolutely continuous functions

$$\{(f(x), g(x)); 0 \leq x \leq 1\}$$

for which $f(0) = g(0) = 0$ and

$$\int_0^1 (\dot{f}^2(x) + \dot{g}^2(x))dx \leq 1.$$

Let $W_1(\cdot)$ and $W_2(\cdot)$ be two independent standard Wiener processes starting from $0$, and let $\eta_2(0,\cdot)$ be the local time process at zero of $W_2(\cdot)$. Then the net of random vectors

$$\left(\frac{W_1(\eta_2(0,xt))}{2^{3/4}t^{1/4}(\log\log t)^{3/4}}, \frac{W_2(xt)}{(2t\log\log t)^{1/2}}; 0 \leq x \leq 1\right)_{t \geq 3}$$

as $t \to \infty$, is almost surely relatively compact in the space $C([0,1], \mathbb{R}^2)$ and its limit points is the set of functions

$$\mathcal{S}^{(2)} := \Big\{(f(h(x)), g(x)): (f,g) \in \mathcal{S}^2, h \in \mathcal{S}_M,$$
$$\int_0^1 (\dot{f}^2(x) + \dot{g}^2(x) + \dot{h}^2(x))\, dx \leq 1, g(x)\dot{h}(x) = 0 \text{ a.e.}\Big\}$$
$$= \Big\{(k(x), g(x)): k(0) = g(0) = 0,\ k, g \in \dot{C}([0,1], \mathbb{R})$$
$$\int_0^1 (|3^{3/4}2^{-1/2}\dot{k}(x)|^{4/3} + \dot{g}^2(x))\, dx \leq 1,\ \dot{k}(x)g(x) = 0 \text{ a.e.}\Big\},$$

where $\dot{C}([0,1], \mathbb{R})$ stands for the space of absolutely continuous functions in $C([0,1], \mathbb{R})$. Consequently : *For the random walk* $\{\mathbf{C}(n) = (C_1(n), C_2(n));\ n = 1, 2, \ldots\}$ *on the 2-dimensional comb lattice* $\mathbb{C}^2$ *we have that*

- *the sequence of random vector-valued functions*

$$\left(\frac{C_1(xn)}{2^{3/4}n^{1/4}(\log\log n)^{3/4}}, \frac{C_2(xn)}{(2n\log\log n)^{1/2}};\ 0 \leq x \leq 1\right)_{n \geq 3}$$

*is almost surely relatively compact in the space* $C([0,1], \mathbb{R}^2)$ *and its limit points is the set of functions* $\mathcal{S}^{(2)}$.

- $\limsup_{n \to \infty} \dfrac{C_1(n)}{n^{1/4}(\log\log n)^{3/4}} = \dfrac{2^{5/4}}{3^{3/4}}$ a.s., *and* $\limsup_{n \to \infty} \dfrac{C_2(n)}{(2n\log\log n)^{1/2}} = 1$ a.s.

- *The limit points of*

$$\left(\frac{C_1(n)}{n^{1/4}(\log\log n)^{3/4}}, \frac{C_2(n)}{(2n\log\log n)^{1/2}}\right)_{n \geq 3}$$

*is not*

$$R = \left[-\frac{2^{5/4}}{3^{3/4}}, \frac{2^{5/4}}{3^{3/4}}\right] \times [-1, 1]$$

*but the domain*

$$D = \{(u, v):\ k(1) = u,\ g(1) = v,\ (k(\cdot), g(\cdot)) \in \mathcal{S}^{(2)}\}.$$

For the many consequences of this result see [R187].

The concept of relative compactness goes back to Strassen [43]. A nice description of this concept is in the *Strong Approximation in Probability and Statistics* book of Pál and Miklós (page 37).

In [R193] we investigated the local time of our random walk on the comb lattice. We define the local time $\Xi(\cdot, \cdot)$ of the random walk $\{\mathbf{C}(n);\ n = 0, 1, \ldots\}$ on the 2-dimensional comb lattice $\mathbb{C}^2$ by

$$\Xi(\mathbf{x}, n) := \#\{0 < k \leq n : \mathbf{C}(k) = \mathbf{x}\}, \tag{3.1}$$

where $\mathbf{x} \in \mathbb{Z}^2$, $n = 1, 2, \ldots$ We proved that:

*On an appropriate probability space for the random walk $\mathbf{C}(n)$ one can construct two independent standard Wiener processes $\{W_1(t); t \geq 0\}$, $\{W_2(t); t \geq 0\}$ with their respective local time processes $\eta_1(t)$ and $\eta_2(t)$ such that for the local time process of the comb walk we have for any $\delta > 0$*

- $$\sup_{x \in \mathbb{Z}} |\Xi((x,0), n) - 2\eta_1(x, \eta_2(0, n))| = O(n^{1/8+\delta}) \quad a.s.,$$

*for any $0 < \varepsilon < 1/4$*

- $$\max_{|x| \leq n^{1/4-\varepsilon}} |\Xi((x,0), n) - \Xi((0,0), n)| = O(n^{1/4-\delta}) \quad a.s.$$

- $$\max_{0 < |y| \leq n^{1/4-\varepsilon}} \max_{|x| \leq n^{1/4-\varepsilon}} |\Xi((x,y), n) - \frac{1}{2}\Xi((0,0), n)| = O(n^{1/4-\delta}) \quad a.s.,$$

*for any $0 < \delta < \varepsilon/2$, where the maximum is taken on the integers.*

We won't list here the many consequences of these results, but would like to mention that the basic ingredient of the strong approximation result is the earlier mentioned result of Pál; the simultaneous approximation of the simple symmetric walk and its local time by their Brownian counterparts from [R77]. Over the years several people generalized this result from simple symmetric walk to more general ones. Just to mention a few, one can consult the paper of Pál with Endre [R82], Borodin [7], Bass and Khoshnevisan [3], [4]. Consequently, under some natural conditions, we were able to extend our results to the following more general type of random walks on the comb lattice: let $\mathbf{C}(0) = \mathbf{0} = (0,0)$, and

$$\mathbf{P}(\mathbf{C}(n+1) = (x, y+k) \mid \mathbf{C}(n) = (x, y)) = p_2(k), \quad (x, y, k) \in \mathbb{Z}^3, \ y \neq 0,$$

$$\mathbf{P}(\mathbf{C}(n+1) = (x, k) \mid \mathbf{C}(n) = (x, 0)) = \frac{1}{2}p_2(k), \quad (x, k) \in \mathbb{Z}^2,$$

and

$$\mathbf{P}(\mathbf{C}(n+1) = (x+k, 0) \mid \mathbf{C}(n) = (x, 0)) = \frac{1}{2}p_1(k).$$

**The general anisotropic model**

As to the general model all the relevant literature is listed in [R195]. One of the most important work was done by Heyde [23], who proved an almost sure approximation for $C_2(\cdot)$ under the condition (3.2) below:

$$n^{-1}\sum_{j=1}^{n} p_j^{-1} = 2\gamma + o(n^{-\eta}), \qquad n^{-1}\sum_{j=1}^{n} p_{-j}^{-1} = 2\gamma + o(n^{-\eta}) \qquad (3.2)$$

as $n \to \infty$ for some constants $\gamma$, $1 < \gamma < \infty$ and $1/2 < \eta < \infty$.

Heyde *et al.* [26] treated the case when conditions similar to (3.2) are assumed but $\gamma$ can be different for the two parts of (3.2) and obtained almost sure convergence to the so-called oscillating Brownian motion. In [R193] we gave a simultaneous approximation for both coordinates of this anisotropic walk:

*Under the condition (3.2), on an appropriate probability space for the random walk $\{\mathbf{C}(N) = (C_1(N), C_2(N)); N = 0, 1, 2, \ldots\}$ one can construct two independent standard Wiener processes $\{W_1(t); t \geq 0\}$, $\{W_2(t); t \geq 0\}$ so that, as $N \to \infty$, we have with any $\varepsilon > 0$*

$$\left|C_1(N) - W_1\left(\frac{\gamma-1}{\gamma}N\right)\right| + \left|C_2(N) - W_2\left(\frac{1}{\gamma}N\right)\right| = O(N^{5/8 - \eta/4 + \varepsilon}) \quad a.s.$$

The anisotropic walk is called periodic if $p_j = p_{j+L}$ for each $j \in \mathbb{Z}$, where $L \geq 1$ is a positive integer. In that case we proved that the rata in the above theorem can be replaced by $O(N^{1/4+\varepsilon})$ and the value of $\gamma$ is $\frac{\sum_{j=0}^{L-1} p_j^{-1}}{2L}$.

Here we present the Strassen theorem, which we proved for the anisotropic walk in [R195]. To begin with we recall some definition. Define the continuous time process $\{\mathbf{C}(u), u \geq 0\}$ by linear interpolation of $\mathbf{C}(N)$. The space $C([0,1], \mathbb{R}^2)$ is the set of continuous functions defined on $[0,1]$ with values in $\mathbb{R}^2$.

*Under condition (3.2) above for the random walk $\mathbf{C}(\cdot)$ we have*

- (i) *the sequence of random vector-valued functions*

$$\left(\sqrt{\frac{\gamma}{\gamma-1}}\frac{C_1(xN)}{(2N \log \log N)^{1/2}}, \sqrt{\gamma}\frac{C_2(xN)}{(2N \log \log N)^{1/2}}, 0 \leq x \leq 1\right)_{N \geq 3}$$

*is almost surely relatively compact in the space $C([0,1], \mathbb{R}^2)$ and its limit points is the set of functions $\mathcal{S}^2$.*

- (ii) *In particular, the vector sequence*

$$\left(\frac{C_1(N)}{(2N \log \log N)^{1/2}}, \frac{C_2(N)}{(2N \log \log N)^{1/2}}\right)_{N \geq 3}$$

*is almost surely relatively compact in the rectangle*

$$\left[-\frac{\sqrt{\gamma-1}}{\sqrt{\gamma}}, \frac{\sqrt{\gamma-1}}{\sqrt{\gamma}}\right] \times \left[-\frac{1}{\sqrt{\gamma}}, \frac{1}{\sqrt{\gamma}}\right]$$

and the set of its limit points is the ellipse

$$\left\{(x, y) : \frac{\gamma}{\gamma - 1} x^2 + \gamma y^2 \leq 1\right\}.$$

- (iii) *Moreover,*

$$\limsup_{N \to \infty} \frac{C_1(N)}{\sqrt{2N \log \log N}} = \frac{\sqrt{\gamma - 1}}{\sqrt{\gamma}} \quad \text{and} \quad \limsup_{N \to \infty} \frac{C_2(N)}{\sqrt{2N \log \log N}} = \frac{1}{\sqrt{\gamma}} \quad a.s.$$

- (iv)

$$\liminf_{N \to \infty} \left(\frac{\log \log N}{N}\right)^{1/2} \max_{1 \leq k \leq N} |C_1(k)| = \frac{\pi \sqrt{\gamma - 1}}{\sqrt{8\gamma}} \quad a.s.$$

$$\liminf_{N \to \infty} \left(\frac{\log \log N}{N}\right)^{1/2} \max_{1 \leq k \leq N} |C_2(k)| = \frac{\pi}{\sqrt{8\gamma}} \quad a.s.$$

In our last joint work, [R200] we explored a randomized version of this model.

In [R194] we studied the local time of the our anisotropic walk as well. Our result were restricted to the above mentioned periodic case. Here we just mention the result without detailing the history and the consequences of them:

*Under the condition (3.2) for the periodic anisotropic walk we have*

$$\mathbf{P}(\mathbf{C}(2N) = (0, 0)) \sim \frac{1}{4\pi N p_0 \sqrt{\gamma - 1}},$$

*for any $L \geq 1$, where* $\gamma = \frac{\sum_{j=0}^{L-1} p_j^{-1}}{2L}$.

Another important issue is the question of recurrence. This was investigated in [R196]. Based on a result of Nash-Williams [36] it was not hard to show that:

*The anisotropic walk is recurrent if*

$$\sum_{k=0}^{\infty} \left(\sum_{j=-k}^{k} \frac{1}{p_j}\right)^{-1} = \infty.$$

On the other hand we proved the following:

*Assume that*

$$\sum_{j=-k}^{k} \frac{1}{p_j} = Ck^{1+A} + O(k^{1+A-\delta}) \quad \text{as } k \to \infty$$

*for some $C > 0, A > 0$ and $0 < \delta \leq 1$. Then the anisotropic random walk is transient.*

**Half Plane Half Comb model**

Beside the special case of the comb model we investigated the case, when we combine the simple symmetric random walk with a random walk on a comb, namely when $p_j = 1/4$, $j = 0, 1, 2, \ldots$ and $p_j = 1/2$, $j = -1, -2, \ldots$, i.e. we have a square lattice on the upper half-plane, and a comb structure on the lower half plane. We call this the Half Plane Half Comb (HPHC) model. For the second component of the HPHC walk a theorem of Heyde et al. [26] gives in this particular case the following strong limit theorem:

*On an appropriate probability space one can construct a sequence $C_2^{(N)}(\cdot)$ and a process $Y(\cdot)$ such that*

$$\lim_{N \to \infty} \sup_{0 \leq t \leq M} \left| \frac{C_2^{(N)}([Nt])}{\sqrt{N}} - Y(t) \right| = 0 \quad a.s.,$$

*where $Y(\cdot)$ is an oscillating Brownian motion (Wiener process) and $M > 0$ is arbitrary.*

In [R193] we gave strong approximation of both components of the random walk $\mathbf{C}(\cdot)$ by certain time-changed Wiener processes (Brownian motions) with rates of convergence. In order to formulate it we need some definitions. Assume that we have two independent standard Wiener processes $W_1(t), W_2(t)$, $t \geq 0$, and let

$$\alpha_2(t) := \int_0^t I\{W_2(s) \geq 0\} \, ds,$$

i.e. the time spent by $W_2(.)$ on the non-negative side during the interval $[0, t]$. Being the process $\gamma_2(t) := \alpha_2(t) + t$ strictly increasing, we can define its inverse: $\beta_2(t) := (\gamma_2(t))^{-1}$. These processes $\alpha_2(t)$, $\beta_2(t)$ and $\gamma_2(t)$ are defined in terms of $W_2(t)$ so they are independent from $W_1(t)$. Here is our theorem:

*On an appropriate probability space for the HPHC random walk $\{\mathbf{C}(N) = (C_1(N), C_2(N)); N = 0, 1, 2, \ldots\}$ with $p_j = 1/4$, $j = 0, 1, 2, \ldots$, $p_j = 1/2$, $j = -1, -2, \ldots$ one can construct two independent standard Wiener processes $\{W_1(t); t \geq 0\}$, $\{W_2(t); t \geq 0\}$ such that, as $N \to \infty$, we have with any $\varepsilon > 0$*

$$|C_1(N) - W_1(N - \beta_2(N))| + |C_2(N) - W_2((\beta_2(N))| = O(N^{3/8+\varepsilon}) \quad a.s.$$

Here the process $W_2(\beta_2(t))$ is identical with $Y(t)$ above i.e. an oscillating Brownian motion. It is a diffusion with speed measure (see [26])

$$m(dy) = \begin{cases} 4\,dy & \text{for} \quad y \geq 0, \\ 2\,dy & \text{for} \quad y < 0. \end{cases}$$

Here again it would be too long to list the consequences of this theorem but let us recall the following LIL type result:

- 
$$\limsup_{t \to \infty} \frac{W_1(t - \beta(t))}{\sqrt{t \log \log t}} = \limsup_{N \to \infty} \frac{C_1(N)}{\sqrt{N \log \log N}} = 1 \quad a.s.,$$

- $$\liminf_{t\to\infty} \frac{W_1(t-\beta(t))}{\sqrt{t\log\log t}} = \liminf_{N\to\infty} \frac{C_1(N)}{\sqrt{N\log\log N}} = -1 \quad a.s.,$$

- $$\limsup_{t\to\infty} \frac{W_2(\beta(t))}{\sqrt{t\log\log t}} = \limsup_{N\to\infty} \frac{C_2(N)}{\sqrt{N\log\log N}} = 1 \quad a.s.,$$

- $$\liminf_{t\to\infty} \frac{W_2(\beta(t))}{\sqrt{t\log\log t}} = \liminf_{N\to\infty} \frac{C_2(N)}{\sqrt{N\log\log N}} = -\sqrt{2} \quad a.s.$$

## 4 Our joint work on spiders

Itô and McKean ([27], Section 4.2, Problem 1) in 1965 introduced a simple but intriguing diffusion process that they called skew Brownian motion, that was extended by Walsh [46] in 1978. Walsh introduced it as a Brownian motion with excursions around zero in random directions on the plane. The random directions are values of a random variable in $[0, 2\pi)$ that are independent for different excursions with a constant value during each excursion. This "definition" can be made precise as, e.g., in Barlow, Pitman and Yor [9]. This motion is now called Walsh's Brownian motion. Following Barlow et al. [10] (see also [16]) one can consider a version of Walsh's Brownian motion which lives on $N$ semi-axis on the plane that are joined at the origin, the so-called Brownian spider. Loosely speaking, this motion performs a regular Brownian motion on each one of the semi-axis and when it arrives to the origin, it continues its motion on any of the $N$ semi-axis with a given probability. Thus, one can construct the Brownian spider by independently putting the excursions around zero of a standard Brownian motion on the $j$-th leg of the spider with probability $p_j$, $j = 1, 2, \ldots, N$ with $\sum_{j=1}^{N} p_j = 1$. A natural discrete counterpart of this motion is a random walk on a spider, i.e., replacing the Brownian motions with simple symmetric random walks on the legs. Hajri [21] studied this issue and proved weak convergence to the Brownian spider.

In the last edition of his *Random Walk in a Random and Non-Random Environment* Pál discussed the corresponding random walk model, which we generalized and extended in [R197] and [R199]. The random walk $\mathbf{S}_n$, $n = 1, 2\ldots$, on the spider $\mathbf{SP}(N)$ (a collection of $N$ half line joining at the origin) can be constructed from a simple symmetric walk $S(n)$, $n = 0, 1, \ldots$ on the line as follows. Consider the absolute value $|S(n)|$, $n = 1, 2, \ldots$, that consists of infinitely many excursions around zero, denoted by $G_1, G_2, \ldots$ Put these excursions, independently of each other, on leg $j$ of the spider with probability $p_j$, $j = 1, 2, \ldots, N$. The first $n$ steps of the spider walk $\mathbf{S}(.)$ is what we obtain this way from the first $n$ steps of the random walk $S(\cdot)$. We denote the Brownian spider on $\mathbf{SP}(N)$, as described above by $\mathbf{B}(t)$, $t \geq 0$, that also starts from the body of the spider, i.e., $\mathbf{B}(0) = 0$. In [R197] we proved a strong approximation theorem between $\mathbf{S}_n$ and $\mathbf{B}(n)$ with the rate $O((n \log \log n)^{1/4}(\log n)^{1/2})$. We established the $n$ step transition probabilities of $\mathbf{S}_n$, and the corresponding transition probabilities for Brownian spider. Then we considered the the following question: How high does the walker go up on a particular leg of the spider. Let $H(j, n)$ denote the highest point reached by the random walk on leg $j$ of the spider in $n$ steps. Formally, let

$$H_M(n) = \max_{1 \leq j \leq N} H(j, n), \qquad H_m(n) = \min_{1 \leq j \leq N} H(j, n).$$

Clearly for $H_M(n)$ we have the LIL and the other LIL of Chung namely

$$\limsup_{n \to \infty} \frac{H_M(n)}{\sqrt{2n \log \log n}} = 1, \quad \text{and} \quad \liminf_{n \to \infty} \left(\frac{\log \log n}{n}\right)^{1/2} H_M(n) = \frac{\pi}{\sqrt{8}} \quad a.s.$$

However the behavior of $H_m(n)$ is much more interesting. We proved that:
*independently of the actual values of $p_1, p_2, \ldots p_N$*

$$\limsup_{n \to \infty} \frac{H_m(n)}{\sqrt{2n \log \log n}} = \frac{1}{2N - 1} \quad a.s.$$

*Moreover let $g(t)$, $t \geq 1$, be a nonincreasing function. Then*

$$\liminf_{n \to \infty} \frac{H_m(n)}{n^{1/2} g(n)} = 0 \quad \text{or} \quad \infty \tag{4.1}$$

*according as $\int_1^\infty g(t)\, dt/t$ diverges or converges.*
For the corresponding result for the Brownian spider, a Strassen theorem and many consequences we refer to [R198]. We would like to mention one more issue here, namely the question of increasing number of legs. For this topic we suppose that each of the $N$ legs are selected with the same probability $1/N$. Denote by $M(n, L)$ the event that the walker on $\mathbf{SP}(N)$ in $n$ steps climbs up at least $L$ high on each legs. Then we proved that:
*for any integer $L \leq \frac{N}{\log N}$, for the $\mathbf{SP}(N)$ we have*

$$\lim_{N \to \infty} \mathbf{P}(M([(cLN \log N)^2], L)) = \left(\frac{2}{\pi}\right)^{1/2} \int_{1/c}^\infty e^{-u^2/2}\, du = P\left(|Z| > \frac{1}{c}\right). \tag{4.2}$$

*where $Z$ is a standard normal variable.* In words, this theorem above gives the limiting probability of the event that, as $N \to \infty$, in $[(cLN \log N)^2]$ steps the walker arrives at least to height $L$ on each of the $N$ legs at least once. In [R199] we investigated the issues of the local time and occupation times on $\mathbf{SP}(N)$ of the spiderwalk and Brownian spider. Let us mention just one of the many results. Denote by $T(j, n)$ the occupation time of leg $j$ in the first $n$ steps of the spiderwalk. Furthermore let

$$T_M(n) = \max_{1 \leq j \leq N} T(j, n), \qquad T_m(n) = \min_{1 \leq j \leq N} T(j, n).$$

We proved that:

*independently of the actual values of $p_1, p_2, \ldots p_N$*

$$\liminf_{n\to\infty} \frac{T_M(n)}{n} = \limsup_{n\to\infty} \frac{T_m(n)}{n} = \frac{1}{N} \quad a.s. \tag{4.3}$$

Again the corresponding result is true for the Brownian spider as well.

## 5 Our work about the distance between random walkers on some graphs

According to the famous paper of Pólya [38] from 1921 two independent random walkers meet infinitely often on $\mathbb{Z}^2$ with probability one. To avoid confusion meeting at the same time at the same place nowadays is called collision, so we will use this term. On $\mathbb{Z}$ two walkers not only collide infinitely many times, but they collide even in the origin infinitely many times. In their landmark paper Dvoretzky and Erdős [15] in 1950 recall the celebrated Pólya result and, among their final remarks, they mention without proof that on $\mathbb{Z}$ three independent random walkers collide (all three together) infinitely often with probability one. A short elegant proof was given for this statement in Barlow, Peres and Sousi [8] in 2012. However four walkers in $\mathbb{Z}$, or three walkers in $\mathbb{Z}^2$, will only collide finitely many times with probability one. Khrishnapur and Peres [34] in 2004 studied this problem on the comb lattice. They proved, that even though the comb is recurrent, two independent random walkers on the comb lattice collide only finitely often with probability 1. Our question was in our paper [R198], that in case the walkers collide only finitely often, then how does their distance grow as a function of time. We established upper class results for the distance of two or more walkers and lower class results for the distance of four or more walkers on $\mathbb{Z}$. Similarly, for the distance in $\mathbb{Z}^2$, we gave upper class results for two or more walkers and lower class results for three or more walkers. Finally, we investigated the distance of two or more walkers on the comb lattice. Here we will only mention our results for the two-dimensional comb lattice. Let $u$ and $v$ be two distinct vertices on $\mathbb{C}^2$. We define their distance as the minimal number of steps on the comb lattice the walker has to take to arrive from $u$ to $v$. Or, equivalently, the distance $d(u,v)$ is the length of the shortest path from vertex $u$ to vertex $v$ in $\mathbb{C}^2$. Formally let $S(n)$ denote the position of the random walker on $\mathbb{C}^2$ at time $n$, then

$$p(u,v) := P(S(n+1) = v | S(n) = u) = \frac{1}{deg(u)} \tag{5.1}$$

and

$$d(u,v) := \min\{k > 0 : P(S(n+k) = v | S(n) = u) > 0\}. \tag{5.2}$$

For $K$ independent walks $\{\mathbf{C}^{(i)}(n) = (C_1^{(i)}(n), C_2^{(i)}(n)) \quad i = 1, 2, ..., K\}$ we investigated

$$D_K^{\mathbb{C}^2}(n) = \max_{i \neq j,\ i,j \leq K} d(\mathbf{C}^{(i)}(n), \mathbf{C}^{(j)}(n)),$$

the maximal distance between the $K$ walkers at time $n$, where $d(x,y)$ was defined above. We proved that:

- *For the distance of $K$ walkers on the comb we have*

$$\limsup_{n\to\infty} \frac{D_K^{\mathbb{C}^2}(n)}{2\sqrt{n\log\log n}} = 1 \quad \text{a.s.} \tag{5.3}$$

- *For $K=2$ and any $\varepsilon > 0$*

$$P\left(D_2^{\mathbb{C}^2}(n) \leq (1+\varepsilon)\frac{2^{9/4}}{3^{3/4}}n^{1/4}(\log\log n)^{3/4} \text{ i.o.}\right) = 1. \tag{5.4}$$

- *For every $\varepsilon > 0$ for n big enough*

$$D_2^{\mathbb{C}^2}(n) > n^{1/4-\varepsilon} \quad \text{a.s.}$$

- *Let $a(n)$ be a nonincreasing nonnegative function. Then, for $K \geq 3$*

$$\sqrt{n}a(n) \in \text{LLC}(D_K^{\mathbb{C}^2}(n)) \tag{5.5}$$

*if and only if*

$$\sum_{n=1}^{\infty}(a(2^n))^{K-2} < \infty. \tag{5.6}$$

# List of Publications of Pál Révész

October 24, 2023

**Books**

[1] *The laws of large numbers*, Akadémiai Kiadó and Academic Press, 1967.

[2] *Strong approximations in probability and statistics*, Akadémiai Kiadó and Academic Press, 1981. (M. Csörgő)

[3] *Random walk in random and non-random environments*, World Scientific Publishing Co., 1990.

[4] *Random walks of infinitely many particles*, World Scientific Publishing Co., 1994.

[5] *Random walk in random and non-random environments*, Second edition, World Scientific Publishing Co., 2005.

[6] *Random walk in random and non-random environments*, Third edition, World Scientific Publishing Co., 2013.

**Papers**

[1a] A Borsuk-féle feldarabolási problémához, *Mat. Lapok*, **7** (1956), 108–111. (A. Heppes)

[1b] Zum Borsukschen Zerteilungsproblem, *Acta Math. Acad. Sci. Hungar.*, **7** (1956), 159–162. (A. Heppes)

[2] A latin négyzet és az ortogonális latin négyzet-pár fogalmának egy új általánosítása és ennek felhasználása kísérletek tervezésére, *Magyar Tud. Akad. Mat. Kutató Int. Közl.*, **1** (1956), 379–390. (A. Heppes)